\documentclass[12pt,a4paper,twoside,final,notitlepage, leqno]{article}
\usepackage[english]{babel}
\usepackage[T1]{fontenc}  
\usepackage{epsfig, graphicx, amssymb}
\usepackage{amsmath,amsthm,epsfig,amsfonts,bbm}  
\newcommand{\monitem}{ \smallskip \noindent $\bullet$ \quad  } 
\newcommand{\moneq}{\vspace*{-6pt} \begin{equation} \displaystyle } 
\newcommand{\moneqstar}{\vspace*{-6pt} \begin{equation*} \displaystyle } 
\newcommand{\monendstar}{\vspace*{-6pt} \end{equation*}   }
\newcommand{\monend}{\vspace*{-6pt} \end{equation}   }
\def\section*#1{}
\usepackage{fancyhdr}
\renewcommand{\headrulewidth}{0pt}
\begin{document} 

\fancypagestyle{plain}{ \fancyfoot{} \renewcommand{\footrulewidth}{0pt}}
\fancypagestyle{plain}{ \fancyhead{} \renewcommand{\headrulewidth}{0pt}}

~

% \vskip -.3 cm   

%%%   \bigskip   \bigskip   \bigskip  
 
\centerline {\bf \LARGE Lattice Boltzmann model approximated  }

 \bigskip 

\centerline {\bf \LARGE 
with finite difference expressions }

%%%%%%%%%%%%%%%%%%%%%%%%%%%%%%%%%%%%%%    03 nov 2014   
\bigskip \bigskip %%%  \bigskip \bigskip 

\centerline { \large    Fran\c{c}ois Dubois$^{a,b}$, Pierre Lallemand$^{c}$, }
\bigskip 
\centerline { \large Christian Obrecht$^{d}$ and Mohamed Mahdi Tekitek$^{e}$ }   

\bigskip  \bigskip \bigskip 

\centerline { \it  \small   
$^a$    Conservatoire National des Arts et M\'etiers,} 

\centerline { \it  \small 
Laboratoire de M\'ecanique des Structures 
et des Syst\`emes Coupl\'es, F-75003,  Paris, France}

\centerline { \it  \small 
$^b$  Department of Mathematics, University  Paris-Sud, B\^at. 425, F-91405 Orsay Cedex, France.}

\centerline { \it  \small  $^c$   Beijing Computational Science Research Center,}

\centerline { \it  \small  
Zhanggguancun Software Park II, Haidian District, Beijing 100094,  China.  } 

\centerline { \it  \small  $^d$  Institut National des Sciences Appliqu\'ees de Lyon,}

\centerline { \it  \small  Centre d'\'Energ\'etique et de Thermique de Lyon (UMR 5008), }

\centerline { \it  \small  Campus La Doua - LyonTech, 69621 Villeurbanne Cedex, France.}

\centerline { \it  \small  $^e$  Department of Mathematics, Faculty of Sciences,} 

\centerline { \it  \small  
University of Tunis El Manar, 2092, Tunis, Tunisia}

\bigskip   
\centerline { \it  \small francois.dubois@math.u-psud.fr, christian.obrecht@insa-lyon.fr, }

\centerline { \it  \small pierre.lallemand1@free.fr,  mohamedmahdi.tekitek@fst.rnu.tn}    
\bigskip   

\bigskip %%  \bigskip 

\centerline {16 April 2016 \footnote {\rm  \small $\,$ Contribution  published in 
 {\it Computers and Fluids}, volume  155, pages 3--8, 2017, 
	doi:10.1016/j.compfluid.2016.04.013. 
  {\it Proceedings} of the  ICMMES Conference,  
Beijing Computational Science Research Center, Beijing, China, 
July 20-24, 2015. Edition 23 February 2018. }}

%%%%%%%%%%%%%%%%%%%%%%%%%%%%%%%%%%%%%%%%%%%%%%%%%%%%%%%%%%%%%%%%%%%%%%%%%%
\bigskip  \bigskip 
\noindent  {\bf Abstract. } \qquad 
We show that the asymptotic properties of the link-wise artificial compressibility method
are not compatible with a correct approximation of fluid properties.
We propose to adapt the previous method through a framework suggested 
by the Taylor expansion method and to replace first order terms in the expansion 
by appropriate three or five points finite differences
and to add non linear terms.                        
The ``FD-LBM'' scheme obtained by this method 
is tested in two dimensions for shear wave, Stokes modes and Poiseuille flow. 
The results are compared with the usual 
lattice Boltzmann method in the framework of multiple relaxation times.

 $ $ \\  [2mm]
   %%  [2mm]
   %%  [2mm]

%%%%%%%%%%%%%%%%%%%%%%%%%%%%%%%%%%%%%%%%%%%%%%%%%%%%%%%%%%%%%%%%%%%%%%%%%%%%%%% 
\bigskip \bigskip  \newpage \noindent {\bf \large   Introduction}  
%%%%%%%%%%%%%%%%%%%%%%%%%%%%%%%%%%%%%%%%%%%%%%%%%%%%%%%%%%%%%%%%%%%%%%%%%%%%%%% 

%%%%%%%%%%%%%%%%%%%%%%%%%%%%%%%%%%%%%%%%%%%%%%%%%%%%%%%%%%%%%%%%%%%%%%%%%%%%%%%  
\fancyhead[EC]{\sc{F. Dubois, P. Lallemand, C. Obrecht and MM. Tekitek }} 
\fancyhead[OC]{\sc{Lattice Boltzmann model approximated with finite differences }} 
%%%%%%%%%%%%%%%%%%%%%%%%%%%%%%%%%%%%%%%%%%%%%%%%%%%%%%%%%%%%%%%%%%%%%%%%%%%%%%%  
%%%%%%%%%%%%%%%%%%%%%%%%%%%%%%%%%%%%%%%%%%%%%  jolie numerotation des pages 
\fancyfoot[C]{\oldstylenums{\thepage}}
%%%%%%%%%%%%%%%%%%%%%%%%%%%%%%%%%%%%%%%%%%%%%  fin jolie numerotation des pages 

\noindent 
Lattice Boltzmann models (LBM)  make it possible  to simulate various types of fluid
flows with simple algorithms (see {\it e.g.} \cite{CD98, De02, LL00, WWLL13,
Ye02}). Usually one can observe (and in simple cases, prove) second order
accuracy (see {\it e.g.} \cite{JKL05}).  These features make LBM approaches
increasingly popular for engineering applications besides others. However, unlike
standard simulation methods such as finite differences, lattice Boltzmann
models are required to process more information than the primitive hydrodynamic
variables, which leads to higher memory consumption and larger data throughput
per collocation point.

\noindent 
On modern computers, especially when using massively parallel processors such as
graphics processing units (GPUs), the computational performance of the LBM is
memory-bound, and therefore is directly linked to the size of the stencil
associated to each collocation point. Asinari {\it et al.} \cite{AOCR12, OAKR14,
OA10, OAY11} proposed the link-wise artificial compressibility method (LW-ACM)
in which parts of the LBM algorithm are replaced by expressions deduced from
finite differencing the primitive variables and gave some results that looked
quite encouraging. Compared to standard three-dimensional LBM, the LW-ACM
reduces memory consumption by a factor of 4.75 and increases performance of GPU
implementations by approximately by a factor of 1.8 \cite{OAKR14}.

\noindent 
We present an analysis of some features of the link-wise artificial
compressibility method of Asinari {\it et al.}, showing possible flaws and then
propose alternative finite difference expressions that allow a significant
improvement of the resulting simulations.

%%%%%%%%%%%%%%%%%%%%%%%%%%%%%%%%%%%%%%%%%%%%%%%%%%%%%%%%%%%%%%%%%%%%%%%%%%%%%%%  section 1 
\bigskip \bigskip   \noindent {\bf \large    1) \quad  Definition of the models }   
%%%%%%%%%%%%%%%%%%%%%%%%%%%%%%%%%%%%%%%%%%%%%%%%%%%%%%%%%%%%%%%%%%%%%%%%%%%%%%%%%%%%%%%%%% 
%

\noindent 
For the sake of simplicity, we start from the usual D2Q9 lattice Boltzmann model 
\cite{LL00}
that allows us to simulate  weakly compressible Navier-Stokes flows. 
Using a planar square grid with collocation points located at $x_{ij}=i\ \delta x$,
$y_{ij}=j\ \delta x$, a fluid is represented by 9 real quantities $f_{ij}^n$ at each of
these grid points.
The LBM simulations involve two steps (collision and propagation)
that we describe following d'Humi\`eres \cite {DDH92,DGK}. 
For the collision at each grid point, one makes a linear transformation of the
quantities $f$ to moments $m$ using an orthogonal matrix $\cal M$ which is shown
below together with a physical interpretation~:
\begin{center}
\begin{tabular}{l|rrrrrrrrr|l}
$\rho$ & 1 & 1 & 1 & 1 & 1 & 1 & 1 & 1 & 1 & density\cr
 $J_x$ & 0 & 1 & 0 & 1 & 0 & 1 &-1 &-1 & 1 & mass flux\cr
 $J_y$ & 0 & 0 & 1 & 0 & 1 & 1 & 1 &-1 &-1 & mass flux\cr
 $E$ & -4 & -1 & -1 & -1 & -1 & 2 & 2 & 2 & 2 & energy\cr
 $XX$ & 0 & 1 &-1 & 1 &-1 & 0 & 0 & 0 & 0 & diagonal stress \cr
 $XY$ & 0 & 0 & 0 & 0 & 0 & 1 &-1 & 1 &-1 & off-diagonal stress \cr
 $q_x$ & 0 &-2 & 0 & 2 & 0 & 1 &-1 &-1 & 1 & energy flux\cr
 $q_y$ & 0 & 0 &-2 & 0 & 2 & 1 & 1 &-1 &-1 & energy flux\cr
 $\epsilon$ & 4 &-2 &-2 &-2 &-2 & 1 & 1 & 1 & 1 & square of energy\cr
\end{tabular}
\end{center}

\noindent
Depending on the simulations to be done, we can conserve only the first 
moment to solve thermal-like problems   %%%%%  modif Christian 03 avril 2016
or we can conserve the three first moments to solve fluid problems for two-dimensional space.
%
%the first (to model thermal problem) or the three first moments (to model fluid problem)
%are not affected by collision (they are the conserved moments).
The others (non conserved) are assumed to evolve as
\begin{equation}
m_k^* = m_k+s_k(m_k^{eq}-m_k),
\end{equation}
where $m_k^{eq}$ is an equilibrium value that is a function of the conserved moments
and $s_k$ a relaxation rate. Note that symmetry considerations are useful to
propose expressions for these equilibrium values.

\noindent
The post-collision moments can also be modified by an external force (gravity,
Coriolis, etc..), preferably following the splitting of Strang \cite{De13}~: applying 
half of the perturbation before collision and half after.
Once the new moments are known, applying ${\cal M}^{-1}$ leads to post-collision
$f^{n*}$.
Propagation is simply obtained through
\begin{equation}
f^{n+1}_{i_0j_0}=f^{n*}_{ij},
\end{equation}
where $i_0$ and $j_0$ are indices of the neighboring grid point corresponding to
the elementary velocity used to define the moments $J_x$ and $J_y$. Thus,
once a velocity set has been chosen, the                                      %%%P
``adjustable'' parameters of a LBM model are the expressions of the equilibrium
values of the non-conserved moments and the values of the relaxation rates.

The analysis of a LBM simulation can be done in several ways. The most popular
is a second order analysis based on the Chapman-Enskog development used in the
kinetic theory of classical gases (see {\it e.g.} \cite{DDH92} or   \cite{LL00}).
This allows to compute the kinematic transport
coefficients (diffusivity for just one conserved moment, shear and bulk viscosities
for three conserved moments). 
It also gives first order expressions for the non-conserved moments.     %%%P
More recently it was proposed to obtain equivalent
equations through Taylor's expansions  \cite{ADGL14, Du08, Du09}, 
which allow to study the effect 
of higher order space derivatives in a
much simpler way than does the Chapman-Enskog development (which makes use of non commuting
matrix products). Finally using the dispersion equation allows to study the linear
stability and gives all the information needed to evaluate the properties of a 
simulation model.

%%%%%%%%%%%%%%%%%%%%%%%%%%%%%%%%%%%%%%%%%%%%%%%%%%%%%%   1-a
\bigskip  \bigskip   \noindent {\bf  \large    1-a) \quad  Standard D2Q9 }    
%%%%%%%%%%%%%%%%%%%%%%%%%%%%%%%%%%%%%%%%%%%%%%%%%%%%%%%%%%%%  

\noindent 
The standard  D2Q9~\cite{LL00} model for Navier-Stokes  uses the following parameters
\begin{center}
\begin{tabular}{cclc}
\cr
moment & & equilibrium &  rate \cr
E      &=& $-2 \rho + 3(J_x^2+J_y^2)/\rho$ &$ s_e$ \cr
XX     &=& $(J_x^2-J_y^2)/\rho$ & $s_{xx}$ \cr
XY     &=& $(J_x\ J_y)/\rho $ & $s_{xx}$ \cr
$q_x$  &=& $-J_x$ &  $s_q$ \cr
$q_y$  &=& $-J_y$ &  $s_q$ \cr
$\epsilon$ &=& $\rho -3(J_x^2+J_y^2)/\rho $& $s_\epsilon $ \cr
\cr
\end{tabular}
\end{center}

\noindent 
This leads to the following properties~:
\begin{center}
\begin{tabular}{lccl}
speed of sound  & $c_s$&=& $\sqrt{\frac{1}{3}}$, \cr
kinematic shear viscosity &  $\nu $&=&  $\frac{1}{3}(\frac{1}{s_{xx}}-\frac{1}{2})$,\cr
kinematic bulk viscosity &  $\zeta $&=&  $\frac{1}{3}(\frac{1}{s_e}-\frac{1}{2}).$\cr
\end{tabular}
\end{center}
\noindent 
The non linear terms lead to the correct advection of shear and acoustic waves.
However, in advective acoustics framework where a uniform velocity $V$ is given,
the LBM method computes the deviation from  this given advection. 
A linear analysis show that low amplitude shear
waves with wave vector parallel to $V$ are damped with an effective kinematic
shear viscosity 
$$
\nu_\textit{eff} = \nu\ (1\ - 3\ V^2).
$$
 The correction is significant as $V$ may be as large as 
$ \, 0.2$ that is typically up to $ \, 0.35 \, $ times the sound speed $c_s.$ 
In the absence of a large velocity, one can easily get higher order terms in the
equivalent equations which allows to determine a shear ``hyperviscosity'' from
the attenuation rate of shear waves at order four  in space derivatives. Previous
work \cite{ADGL14, DL09, DL11}
showed which conditions allowed to get an isotropic hyperviscosity (no
angular dependence in the expressions) and the possibility to make it equal to zero.

%%%%%%%%%%%%%%%%%%%%%%%%%%%%%%%%%%%%%%%%%%%%%%%%%%%%%%   1-b
\bigskip  \bigskip   \noindent {\bf  \large 
  1-b) \quad  Link-wise artificial compressibility method }    
%%%%%%%%%%%%%%%%%%%%%%%%%%%%%%%%%%%%%%%%%%%%%%%%%%%%%%%%%%%%%%%%%%%%%%%%%%%%%%%%%%%%  

\noindent
The new proposal of Asinari {\it et al.}  \cite{AOCR12, OAKR14, OA10, OAY11}
uses just the primitive variables: density $\rho$ and velocity $\vec u$. From
these quantities it reconstructs a set of $f^n$ on all grid points of the
computation domain and then lets them evolve with the LBM rules. In its original
formulation, the reconstruction rule is expressed through the equilibrium
distribution $f^{eq}$ which is function of the sole primitive variables. Using
the present notations, it can be written as:
\begin{equation}
f^{n*}_{ij} = f^{eq}(\rho^n_{ij},\vec{u}^{\,n}_{ij})
+\Theta\left(f^{eo}(\rho^n_{i_0j_0},\vec{u}^{\,n}_{i_0j_0})
-f^{eo}(\rho^n_{ij},\vec{u}^{\,n}_{ij})\right) \, ,
\end{equation}
where $f^{eo}(\rho,\vec{u})$ is defined as:
$$f^{eo}(\rho,\vec{u})=\frac{1}{2}\left(f^{eq}(\rho,\vec{u})-f^{eq}(\rho,-\vec{u})\right) \, ,$$
and $\Theta$ as:
$$\Theta=1-\frac{2\nu}{c_s^2}=1-6\nu \, .$$
The properties of the proposed algorithm lies in the reconstitution. The work
of Asinari {\it et al.}  use what can be called ``zeroth-order'' reconstitution
as they just involve the expressions shown in the preceding table.

\noindent 
To analyze it we use a classical Von Neumann stability 
analysis in Fourier space (see \cite{LL00}).
So we proceed in the following way. Starting either from the equations
to be simulated or from the computer code derived from them we prepare a series of
instructions for  a computer algebra system. We then consider a 
grid with the following initial conditions~: a plane                       %%%P
wave of small amplitude and wave vector $k_x,\ k_y$,
uniform density plus possibly a  uniform velocity $V=(V_x,V_y).$ 
This means we take the following initial state~:
$f=f^0 + \delta f$, where $f^0=(f_0,\dots,f_8)$ 
represents  the uniform equilibrium state specified by uniform and steady density 
$\rho$ and velocity $V=(V_x,V_y)$ and $\delta f=(\delta f_0,\dots,\delta f_8)$ is the fluctuation.
We then apply one time step in the Fourier space and linearize the results
in terms of the parameters of the plane wave (amplitude and phase factors).

\noindent 
We define space phase factors $p= e^{i\,  k_x }$ and $q=e^{i\, k_y}$ and time factor
$z=e^{-\Gamma}$ ($i$ is unit imaginary number and $\Gamma$~being the attenuation rate) 
in units such that $\delta x=1$
and the duration of one time step equals to unity.
So the initial conditions in moment space are
$$  
\delta \rho(j,l) \,=\, A\ p^j\ q^l \,, \quad 
\delta J_x(j,l)\,=\,B \, p^j\ q^l  \,, \quad 
\delta J_y(j,l) \,=\, C\ p^j\ q^l  .
$$ 
In consequence we have classical relation of the type 
$$\delta \rho (j+1,l)= e^{i k_x} \, \delta \rho(j,l) = 
p \, \delta \rho(j,l), \quad \delta \rho (j,l+1)= e^{i k_y} \, \delta \rho(j,l)=q \, \delta \rho(j,l),$$
and  analogous relations for two  others fields 
$\delta J_x$ and $\delta J_y.$ We introduce the state vector
$\Phi=\left (A,B,C\right)^{{\rm{t}}}$, after one time step the vector $\Phi$ 
is multiplied by the amplification matrix $H$~:
\begin{equation}
\Phi^{n+1}=H \, \Phi^n.
\label{iteration}
\end{equation}
We note here that the amplification matrix $H$ is determined by the collision step and the advection step.
In particular the coefficients~: $V$, $c_s$, $s_k$, $p$ and  $q$ (see for details the original 
reference \cite{LL00}). 
%after 1 time step, one can write the state 
$\Phi=\left (\delta \rho,\ \delta J_x,\ \delta J_y \right)^{{\rm{t}}}$
%in matrix form
We search the modes associated to the iteration (\ref{iteration}). In that case the vector $\Phi$ is solution
of 
\begin{equation}
z \ \Phi= H \ \Phi \, ,
\end{equation}
from which we get the dispersion equation 
\begin{equation}
F(p,q,z)=\det(H-z\ Id)\, ,
\label{disp}
\end{equation}
where $Id$ is the unit matrix.
Literal expressions of $F$ are then solved to get $z$ by successive approximations in powers
of the wave vector components.  In other terms we search the eigenvalues 
$z$ and eigenvectors $R$ as powers of the wave vector $(k_x,k_y).$
So that the attenuation rate (possibly complex for
propagating waves) is obtained as an expansion in wave vector components
with $\Gamma=\log{z}$.
As the general expressions are quite cumbersome, we only give information on the terms up to
power 2 in wave vector components. In addition we assume that the uniform velocity is
parallel to the wave vector with amplitude $v$  
({\it{i.e.}} $V= (V_x , V_y)$ and $v =|V|$) 
and we apply a rotation of the axis so that
the wave vector is parallel to Ox (rotated axis) with amplitude $k$. 
We replace the spatial phase factor $p$ and $q$ 
by their expansion at second order in in $k$.
We then get the matrix $H(k)$~:
$$
H(k) \,=\, \left( \begin {array} {ccc} 
1-6 \,i\,  v k\nu-\frac{1}{6} k^2(1+3 v^2)  &  -6 \,i\, \nu k-v k^2 & 0  \cr
-\frac{1}{3} \,i\,  k (1+3 v^2)+v (1- 3 v k^2) &  1-2  \,i\,  k v - 3 \nu k^2  & 0 \cr
0 & 0 &  1-  \,i\,  k v - \nu k^2 
\end {array} \right)  . 
$$
Note that no angle appears, so the model is isotropic at order 2 in wave vector. The
previous matrix shows decoupling of one shear mode and two longitudinal modes.

\noindent 
From the roots of the dispersion equation in the case $v=0$, one obtains the kinematic shear
viscosity $\nu$, related to the relaxation rate $s_{xx}$ by
$$ 
\nu=\frac{1}{3} \, \Big(\frac{1}{s_{xx}}-\frac{1}{2} \Big) \,,
$$ 
and the speed of sound and its damping
$$ 
c_s=\sqrt{2 \nu},\quad\ \Gamma_s=\frac{\nu}{2}+\frac{1}{12} \, . 
$$ 
Note that the result for the damping of sound can be interpreted with a kinematic
bulk viscosity independent of the parameters of the model.

\noindent 
 When $v$ is not zero, since the transport coefficients 
can be obtained through a perturbation analysis, 
we shall use the following series expansion in $k$ of the roots \cite{LL00}. 
One can verify that the roots contain a linear  dependence in $v$ 
(term in $ \, i  k v$ linked to linear advection)
and the shear viscosity becomes
$$ 
\nu(v) = \nu - \frac{1}{2}\ v^2 \, .
$$ 
This last result means that if $v>\sqrt{2 \nu}\ =\ c_s$, shear waves grow exponentially
and thus the model is unstable so it is not recommended to use this model
for simulations at fairly large Reynolds number.
Actual simulations allow to verify the previous results (see section 3-a)     %%%P

%%%%%%%%%%%%%%%%%%%%%%%%%%%%%%%%%%%%%%%%%%%%%%%%%%%%%%%%%%%%%%%%%%%%  section 2 
\bigskip  \bigskip   \noindent {\bf  \large    2) \quad  New proposition  }    
%%%%%%%%%%%%%%%%%%%%%%%%%%%%%%%%%%%%%%%%%%%%%%%%%%%%%%%%%%%%%%%%%%%%%%%%%  

\noindent 
We propose to use the same basic idea (reconstruction of the $f$ from primitive
variables~: density, components of the velocity), but with improved formulae.

In the Taylor expansion analysis leading to the equivalent equations  \cite{Du08}, 
it was shown that the 
non conserved moments $m_k$ can be expanded in powers of the size of the elementary step of 
the algorithm. Beyond the order 0, presented above, the second order has been expressed
in terms of $\theta_k$ that involve space derivatives and non linear terms. 
In fact, as described in \cite{Du08}, 
we have the following development of non-equilibrium moments at second order on $\Delta t$~:
\begin{equation}
m_{k}^{*}=m_{k}^{eq}+\Delta t
 \left(\frac{1}{2}-\sigma_{k}\right)\theta_{k}+{\rm{O}}(\Delta t^{2}),
 \quad k \geq 2.
\label{order2}
\end{equation}
where $\sigma_{k}\equiv \left(\frac{1}{s_{k}}-\frac{1}{2}\right)$ and $\theta_k$ 
is the defect of conservation defined by~:
\begin{equation}
\theta_k\equiv \partial_t m_k^{eq} + \Lambda_{k \, \alpha}^{\ell} \partial_{\alpha} m_{\ell}, \quad k > N,
\label{theta}
\end{equation}
where $N$ is the number of the conserved moments and 
$\Lambda_{k \, \alpha}^{\ell}=\sum_j v_j^{\alpha} v_j^{\beta} (M^{-1})_{j k}$, $k
=0\dots 8$, $\alpha=1\dots 2$ and $\beta=1\dots 2$. \\
{\bf{Remark}}
In the case of the $N=3$ ({\it{i.e.}} 3 conserved moment to model fluid-like problems), 
we get the following macroscopic equations~:
$$
\partial_t m_k +  \Lambda_{k \, \alpha}^{\ell} \, \partial_{\alpha} m_{\ell}^{eq} 
- \sigma_{\ell} \, \Delta t \, \Lambda_{k \, \alpha}^{\ell} \, \partial_{\alpha} \theta_{\ell} 
= {\rm{O}} (\Delta t)^2, \quad k=0,1,2.
$$
We note here that for $k=1,2$ at the order one we have a term $\frac{1}{3} \nabla \rho$
 which gives the sound speed $c_s=\frac{1}{\sqrt{3}}.$ At the order two 
(terms  having $\Delta t$ as coefficient) we obtain the viscous terms function 
of $\sigma_{\ell}$. For more details see \cite{Du08}.
 
\noindent As many individual terms are found to play no role in the behavior of the shear and
acoustic modes, we give only the relevant terms of the defect 
of conservation $\theta_k$ (\ref{theta}) for the case where the density is close to 1~:
$$
\left \{
\begin{array}{lcccl}
\theta_3 &\equiv& \theta_E  & = &    (2 + 6\ (v_x^2+v_y^2))\ (\partial_x v_x+\partial_y v_y)
          -2\ (v_x\partial_x\rho+v_y\partial_y\rho),\\
\theta_4 &\equiv& \theta_{XX} & = & \frac{2}{3}(\partial_x v_x-\partial_y v_y)
  -\frac{2}{3} (v_x \partial_x \rho-v_y \partial_y \rho) \\
& & &   &   \qquad - \,  2(v_x(\partial_x v_x^2+\partial_y v_xv_y)
            -v_y(\partial_x v_xv_y+\partial_y v_y^2) ,\\
\theta_5 &\equiv&   \theta_{XY} & = & \frac{1}{3}(\partial_x v_y-\partial_y v_x)
  -\frac{1}{3} (v_x \partial_y \rho+v_y \partial_x \rho) \\
& & & &   \qquad - \, v_x(\partial_x v_xv_y+\partial_y vy^2) 
      -v_y(\partial_x v_x^2+\partial_y v_x v_y) \, . 
\end{array}
\right.
$$

The partial derivatives are then estimated by finite difference.
To sum up, the neighboring $f$ are obtained (see equation (\ref{order2})) using the non-conserved moments~:
$$
\left \{
\begin{array}{lcl}
E & = &  -2 \rho + 3(v_x^2+v_y^2)/\rho +(1-\frac{1}{s_e})\ \theta_E ,\\
XX  & = & (J_x^2-J_y^2)/\rho + (1-\frac{1}{s_{xx}})\ \theta_{XX} ,\\
XY  & = &  (J_x\ J_y)/\rho + (1-\frac{1}{s_{xx}})\ \theta_{XY} ,\\
q_x & = & -J_x, \\
q_y & = & -J_y, \\
\epsilon & = &   \rho - 3(v_x^2+v_y^2)/\rho. \\
\end{array}
\right.
$$
\noindent 
With these expressions, the acoustic waves propagate with speed $1/\sqrt{3}$
(as for standard D2Q9), advection by a mean velocity $V$ is correct and the
viscosities are now~:                                        %%%P
%\begin{center}
%\begin{tabular}{lcl}
%\cr
$$\mbox{shear } =  \frac{1}{3}(\frac{1}{s_{xx}}-\frac{1}{2})\ (1-3\ V^2) \quad \mbox{ and } \quad 
\mbox{ bulk } = \frac{1}{3}(\frac{1}{s_e}-\frac{1}{2})\ (1-3\ V^2), 
$$
%\cr
%\end{tabular}
%\end{center}
%\noindent 
as is known for D2Q9.

\noindent 
For the particular case with $V=0$, one can determine higher order contributions   %%%P
to the damping of the hydrodynamic modes. 
 We first expand the dispersion equation (\ref{disp}),
then we replace spatial phase factors $p$ and $q$ by their 
expansions up to the fourth order in $k$ and solve the resulting expression 
by successive approximation in $k$. This leads to eigenvalues $z_i,\,  i=1..3$ 
and then we get the development of the damping 
coefficient $\Gamma_i=-\log(z_i)$. We interpret one of these roots as 
$$\Gamma_i=\nu_0 k^2+ \nu_2 \, k^4.$$
Which allows to define a $k$ dependent kinematic shear viscosity~:
$$ \nu(k)=\nu_0 + \nu_2 \, k^2, $$%+  \nu_2 \, k^4  - \nu_4 \, k^6 + \dots  $$
We define the coefficient  $\nu_2$ as
 ``hyperviscosity''.
 The expressions
for this hyperviscosity depend on the way space derivatives are estimated using
finite difference.

We have considered three cases.

\smallskip \monitem {\bf Three points stencil} \quad such that 
$$ 
\partial_x \bullet \,\simeq \, \frac{1}{2} \, \big( \bullet (i+1,j)-\bullet (i-1,j) \big)  \,. 
$$ 
Then the shear hyperviscosity is
$$ 
\nu_2=\frac{1}{72}(2\sigma_{xx}-3)(2\sigma_{xx}-1)-\frac{8\sigma_{xx}-3}{36}
(\cos{\phi}^2-\cos{\phi}^4),
$$ 
where $\sigma_{xx}=1/s_{xx}-1/2$ and $\phi$ is the angle between the Ox axis and the
wave vector. This contribution is anisotropic. It becomes larger than the usual
viscous term for $k>0(\sqrt{\sigma_{xx}})$ which will prevent from doing significant
simulations at small viscosity.

\smallskip \monitem  {\bf Five points stencil}  \quad such that 
$$ 
\partial_x \bullet  \,\simeq \, 
\frac{3}{4} \, \big(\bullet (i+1,j)-\bullet (i-1,j) \big) 
\,- \, \frac{1}{8} \, \big(\bullet (i+2,j)-\bullet (i-2,j) \big) \,. 
$$ 
This leads to a shear hyperviscosity
$$ 
\frac{1}{36}\sigma_{xx}(2\ \sigma_{xx}-1)-\frac{\sigma_{xx}}{18}(\cos{\phi}^2-\cos{\phi}^4) .
$$ 
This is still anisotropic but removes the small viscosity limitation.

\smallskip \monitem  {\bf Nine points stencil} \quad based on the D2Q9 geometry, we can use
$$ 
\begin{array}{rcl}
\partial_x \bullet  & \,  \simeq  \,& 
\bullet (i+1,j)-\bullet (i-1,j) \\ & \,-\, & \frac{1}{4} \,
\big[\bullet (i+1,j+1)-\bullet (i-1,j+1)-\bullet (i-1,j-1)+\bullet (i+1,j-1) \big]
\end{array}
$$ 
and similar expression for $\partial_y $.
This leads to the following shear hyperviscosity~:
$$ 
\frac{1}{24\ \sigma_{xx}}(3-2\ \sigma_{xx})(2\ \sigma_{xx}-1)
-\frac{20\ \sigma_{xx}-9}{12 \sigma_{xx}}(\cos{\phi}^2-\cos{\phi}^4)
$$ 
which is still anisotropic and does not solve the limitation indicated for
the three point stencil.

\noindent 
For all three stencils, the full dispersion equation (cubic equation in time factor
$z$) can be obtained numerically for $k$ up to $\pi$ in order to predict the linear
stability.

%%%%%%%%%%%%%%%%%%%%%%%%%%%%%%%%%%%%%%%%%%%%%%%%%%%%%%%%%%%%%%%  section 3 
\bigskip \bigskip   \noindent {\bf \large  3) \quad  Results of some simulations  }   
%%%%%%%%%%%%%%%%%%%%%%%%%%%%%%%%%%%%%%%%%%%%%%%%%%%%%%%%%%%%%%%%%%%%%%%%%%%%%%%%%%%%%%%  

%%%%%%%%%%%%%%%%%%%%%%%%%%%%%%%%%%%%%%%%%%%%%%%%%%%%%%%%%%%%%%%  section 3-a 
\bigskip   \noindent {\bf \large   3-a) $ \,\,\, $ Shear wave  }    
%%%%%%%%%%%%%%%%%%%%%%%%%%%%%%%%%%%%%%%%%%%%%%%%%%%%%%%%%%%%%%%%%%%%%%%%%%%%%%%%%%%%%%%  

\noindent 
Elementary tests have been performed in a square domain (size $N^2$) with periodic
boundary conditions. The initial condition is a shear wave of wave vector $k_x,k_y$
(of modulus $k$)
with in some cases a uniform velocity parallel to the wave vector. 
 In fact we take the following initial conditions~:
$$\left \{ \begin{array}{rcl}
\rho(t=0)&=&1,\\
j_x(t=0)&=&- A(0) \frac{k_y}{k} \cos(k_x \, x +k_x\, y) + \frac{k_x}{k} V,\\
j_y(t=0)&=& A(0) \frac{k_x}{k} \cos(k_x \, x +k_x\, y) + \frac{k_y}{k} V,
\end{array}
\right.$$
The exact solution admits the same algebraic form,    %%% Christian, 3 avril 2016 
except that $A$ is replaced by a function of time $A(t)$ ; then $A=A(0)$. 
At each time
step we measure the correlation function  $A(t)/A(0)$ of the velocity field with its initial
state. For $V=0$,  $A(t)$  decays exponentially, otherwise it is 
$ \,  e^{(-\Gamma t)} \,  \big( \cos{\omega t} \big)  $.

%
%%%%%%%%%%%%%%%%%%%%%%%%%%%%%%%%%     Figure 1        %%%%%%%%%%%%%%%%%%%%%%%%%%%%%% 
\begin{figure}[htp]  \begin{center} 
\label{ACM_original}  
\centerline { \includegraphics[width=.75 \textwidth] {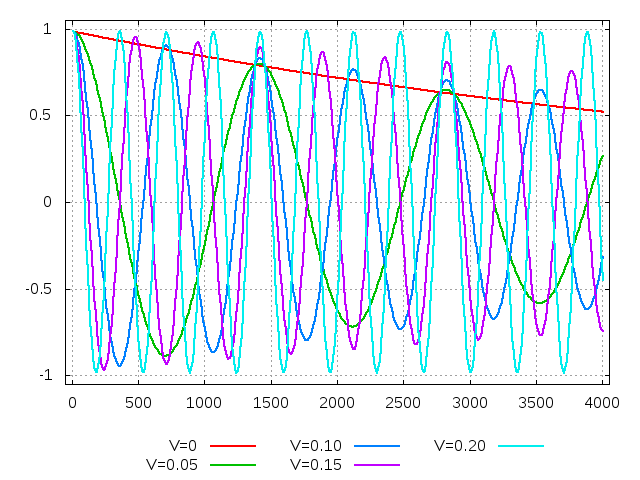}}
\caption{Time evolution of the correlation function  $A(t)/A(0)$ {\it{versus}} discrete time
between  an initial transverse wave (of vector $(3\, k_0 ,2 \, k_0)$
 where $k_0=\frac{2 \pi}{191}$) 
and its later state for five different values of the mean velocity $V$.   %%%P
Square of 191 $\times$ 191 nodes and periodic boundary conditions.       %%%P
When  $V$  grows, the dissipation of the waves is reduced. } 
   \end{center} \end{figure}     
%%%%%%%%%%%%%%%%%%%%%%%%%%%%%%%%%     Figure 1        %%%%%%%%%%%%%%%%%%%%%%%%%%%%%% 
%

\noindent 
We show in Fig.~1 %% \ref{ACM_original} 
the results for the initial ACM model (no $\theta$ in our proposal)
for 5 values of the mean velocity $V$. Clearly the velocity square dependence of   %%%P
the damping is unacceptable.

\noindent 
We then perform a series of measurements at $V=0$ for several values of the wave vector
and compare (Table 1)  the measured relaxation rate $\Gamma$ to the development in terms of
hyperviscosity and the numerical root of the dispersion equation (that which
corresponds to the transverse mode). Fig.~2,~3 and~4 
illustrate the results for the three, five and nine points stencil respectively. These figures
have been obtained for a long wave length kinematic shear viscosity $\nu_0=0.01$.
In the case of the nine point stencil, the model is unstable in the \{1,1\} direction
so no simulation could be performed.

In fact in table \ref{hyper} we study the equivalent hyperviscosity for the 
ACM scheme for different stencils.   %%%%%%   Christian    03 avril 2016 
We show that the hyperviscosity is relatively high and negative for an
angle equal to $45^{\rm o}$ for 
the nine point stencil. This is directly correlated to instability in the \{1,1\} direction. 

%%%%%%%%%%%%%%%%%%%%%%%%%%%%%%%%%%%%%%%%%%%%%%%%%%%%%%%%%%%%%%%%%%%%%%%%%%%%%%%
\begin{table}
\begin{centering}
\begin{tabular}{cccc}
angle & 3-point&  5-point  &  9-point \cr
0.00    &  0.03725 & -0.00103 &   0.03725\cr
 26.60 &   0.02534 & -0.00067  &  0.00079\cr
45.00  &   0.01867 & -0.00047  &-0.0196\cr
 \end{tabular}
 \caption{Numerical study of the hyperviscosity for different stencils
 of ACM scheme {\it{vs}} angle of the wave vector $k.$ 
 All simulations are performed with the same value $s_{xx}=1.85$.}
\end{centering} 
\label{hyper}
\end{table}
%%%%%%%%%%%%%%%%%%%%%%%%%%%%%%%%%%%%%%%%%%%%%%%%%%%%%%%%%%%%%%%%%%%%%%%%%%%%%%%

%
%%%%%%%%%%%%%%%%%%%%%%%%%%%%%%%%%     Figure 2    %%%%%%%%%%%%%%%%%%%%%%%%%%%%%% 
%%%%%% \begin{figure}[htp] 
\begin{figure}[!h] \begin{center} 
\label{stencil3}  %%  
\centerline { \includegraphics[width=.72 \textwidth] {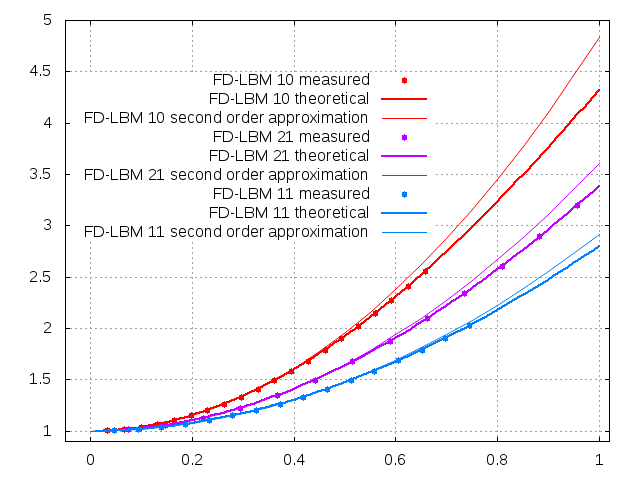} }
\caption{Relative shear viscosity (normalized by  $\nu_0$) for the three point stencil {\it versus} 
 wave vector modulus $k$.
Solid curves from dispersion equation, 
thin solid curves from hyperviscosity, discrete points from actual simulation.    %%%P
Top curves for wave vector along X axis \{1,0\}, middle curves for wave
vector along \{2,1\} direction and lower curves for wave vector along the \{1,1\} direction.}
   \end{center} \end{figure} 
%%%%%%%%%%%%%%%%%%%%%%%%%%%%%%%%%%%%%%%%%%%%%%%%%%%%%%%%%%%%%%%%%%%%%%%%%%%%%%%%%%%%

%
%%%%%%%%%%%%%%%%%%%%%%%%%%%%%%%%%     Figure 3      %%%%%%%%%%%%%%%%%%%%%%%%%%%%%% 
%%%%%% \begin{figure}[htp] 
\begin{figure}[!h] \begin{center} 
\label{stencil5}  
\centerline { \includegraphics[width=.72 \textwidth] {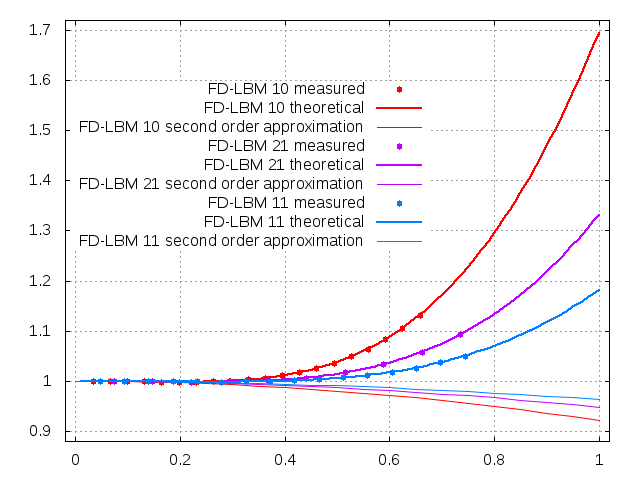}}
\caption{Relative shear viscosity (normalized by  $\nu_0$) for the five  point stencil {\it versus} 
 wave vector modulus $k$.
Solid curves from dispersion equation, 
thin solid curves from hyperviscosity, discrete points from actual simulation.   %%%P
Top curves for wave vector along X axis \{1,0\}, middle curves for wave
vector along \{2,1\} direction and lower curves for wave vector along the \{1,1\} direction.}
\end{center} \end{figure} 
%%%%%%%%%%%%%%%%%%%%%%%%%%%%%%%%%%%%%%%%%%%%%%%%%%%%%%%%%%%%%%%%%%%%%%%%%%%%%%%%%%%%

%%%%%%%%%%%%%%%%%%%%%%%%%%%%%%%%%     Figure 4      %%%%%%%%%%%%%%%%%%%%%%%%%%%%%% 
%%%%%% \begin{figure}[htp] 
\begin{figure}[!h]  \begin{center} 
\label{stencil9}  
\centerline { \includegraphics[width=.73 \textwidth] {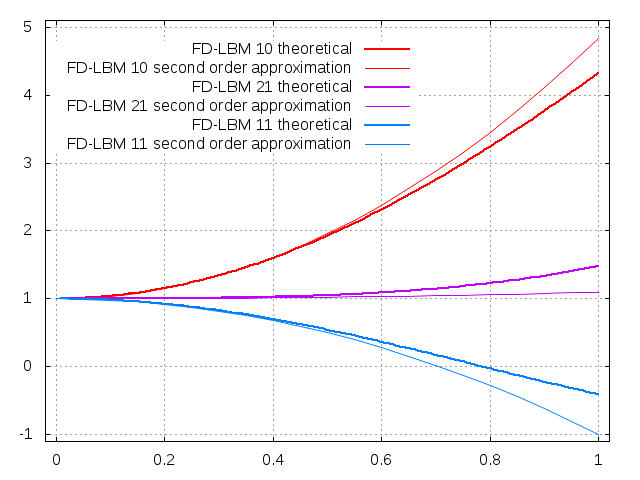}}
\caption{Relative shear viscosity (normalized by  $\nu_0$) for the nine point stencil {\it versus} 
 wave vector modulus $k$.
Solid curves from dispersion equation, thin solid curves from hyperviscosity, squares from %%%P
actual simulation. Top curves for wave vector along X axis, middle curves for wave
vector along \{2,1\} direction and lower curves for wave vector along the \{1,1\} direction.
No experimental data due to linear instability at least  in the \{1,1\} direction. }   %%%P
\end{center} \end{figure} 
%%%%%%%%%%%%%%%%%%%%%%%%%%%%%%%%%%%%%%%%%%%%%%%%%%%%%%%%%%%%%%%%%%%%%%%%%%%%%%%%%%%%

%
%%%%%%%%%%%%%%%%%%%%%%%%%%%%%%%%%     Figure 5        %%%%%%%%%%%%%%%%%%%%%%%%%%%%%% 
%%%%%% \begin{figure}[htp] 
\begin{figure}[!h] \begin{center} 
\label{figd2q9}  
\centerline { \includegraphics[width=.73 \textwidth] {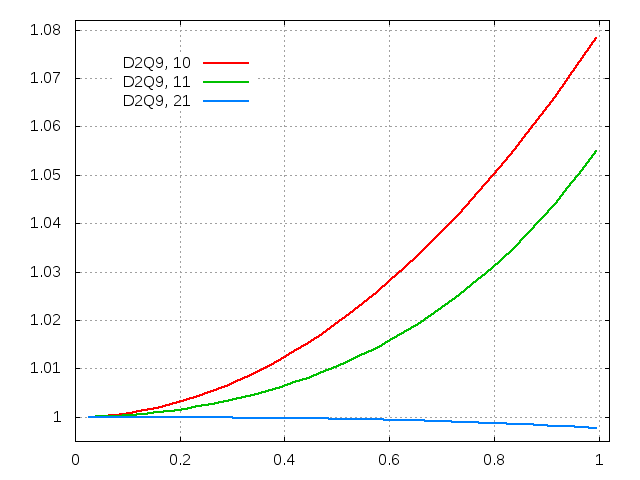}}
\caption{Relative shear viscosity for the D2Q9 lattice Boltzmann model {\it versus} 
 wave vector modulus $k$  for three orientations of the wave vector 
(top   \{1,0\}, middle \{2,1\}, bottom \{1,1\}).
We observe that the maximal error (8 \%) is very much reduced compared to the FD-LBM
scheme presented at Figures 2, 3 and 4.} 
\end{center} \end{figure} 
%%%%%%%%%%%%%%%%%%%%%%%%%%%%%%%%%%%%%%%%%%%%%%%%%%%%%%%%%%%%%%%%%%%%%%%%%%%%%%%%%%%%
% 

%%%%%%%%%%%%%%%%%%%%%%%%%%%%%%%%%%%%%%%%%%%%%%%%%%%%%%%%%%%%%%%  section 3-b 
\bigskip \bigskip  \newpage  \noindent {\bf \large   3-b) $ \,\,\, $ Stokes modes  }    
%%%%%%%%%%%%%%%%%%%%%%%%%%%%%%%%%%%%%%%%%%%%%%%%%%%%%%%%%%%%%%%%%%%%%%%%%%%%%%%%%%%%%%%  

%%%
\noindent 
We give some partial results of simulations of situations less elementary that
simple plane waves. To take solid boundaries into account we propose to consider
the lattice nodes just outside the fluid region and to estimate the state of the
virtual fluid in those points by linear extrapolation using the fact that the 
velocity is 0 on the boundary. As in the scheme of Bouzidi et al.\cite{BFL01}    
stability is obtained by using different expressions depending on the location
of the intersection of the boundary with the link that goes from 
the last fluid point to the first solid point.  %%%P

\noindent 
We then compute the relaxation rate of the Stokes modes inside a circle of radius
$R=29.9$ lattice units. The flow field is obtained from the stream function
\begin{equation}
\psi(r,\theta,t)=e^{-\big(\displaystyle \Gamma\,t\big)}\ \cos{(m\ \theta)}J_n( r/R),
\end{equation}
with singlets for $m=0$ and doublets for $m>0$ and
\begin{equation}
\Gamma= \frac{\nu}{R^2} \, a_l^2 \, ,
\end{equation}
where $a_l$ is a $l^{th}$ zero of the Bessel functions  $J_m(a_l)=0$.
We give in the following table some values of the relative difference between
measured values and the theoretical values for three cases~: present FD-LBM
with the three point stencil, optimized LBM-D2Q9 ($\nu=1/\sqrt{108}$) and 
\begin{equation} \label{sigsig} 
\Big( \frac{1}{s_{xx}}-\frac{1}{2} \Big) \ \Big( \frac{1}{s_{q}}-\frac{1}{2} \Big)
=\frac{1}{6} \, ,
\end{equation}
required to yield an isotropic hyperviscosity), and a non-optimized D2Q9-LBM
(same value of $\nu$, but $s_{q}=1.3$ instead of 0.9282. It is clear that
FD-LBM does not match the accuracy of optimized LBM-D2Q9 (see \cite{DL09}).

\begin{table}
\begin{centering}
\begin{tabular}{cccccc}
$l$ & Bessel &  FD-LBM-3  &  FD-LBM-5 &   BGK  &  LBE-q   \cr
&&&Singlets&&\cr
 1 &  14.68200   &  0.00729 &   0.00003  &  0.00053  &  -0.00010   \cr
 2 &  49.21850   &  0.02191 &  -0.00141  &  0.00179  &  -0.00114   \cr
 3 & 103.49950   &  0.04663 &  -0.00313  &  0.00382  &  -0.00276   \cr
 4 & 177.52080   & 0.07969  &  -0.00400  &  0.00672  &  -0.00489   \cr
 5 & 271.28171   &  0.12335 &  -0.00358  &   0.01071 &  -0.00752   \cr 
 6 & 384.78189   & 0.17778  &  -0.00099  &  0.01623  &  -0.01053   \cr 
&&&Doublets&&\cr

 1 &  26.37460  & 0.01324 &  -0.00090  &  0.00106  &  -0.00042  \cr
 2 &  40.70650  & 0.02078 &  -0.00103  &  0.00164  &  -0.00087  \cr
 3 &  57.58290  & 0.02959 &  -0.00147  &  0.00236  &  -0.00133  \cr
 4 &  76.93890  & 0.03966 &  -0.00186  &  0.00323  &  -0.00183  \cr
 5 &  98.72630  & 0.05060 &  -0.00231  &  0.00424  &  -0.00236  \cr
 6 & 122.90760  & 0.06241 &  -0.00254  &  0.00538  &  -0.00293  \cr
 7 & 149.45290  & 0.07545 &  -0.00275  &  0.00667  &  -0.00354  \cr
 8 & 178.33730  & 0.08948 &  -0.00267  &  0.00809  &  -0.00419  \cr
 9 & 209.54010  & 0.10418 &  -0.00230  &  0.00965  &  -0.00488  \cr
10 & 243.04340  & 0.12003 &  -0.00175  &  0.01138  &  -0.00563  \cr
11 & 278.83160  & 0.13682 &  -0.00099  &  0.01328  &  -0.00643  \cr
 \end{tabular}
 \caption{Numerical study of the Stokes modes in a disk. 
 All simulations are performed with the same value $\nu_0=0.05$.
The second  column gives the theoretical values and the other the error 
between the numerical scheme and the theoretical value. 
The third column uses the present FD-LBM scheme with a three point
stencil for the evaluation of the gradients, 
the fourth column  the present FD-LBM scheme with a five point
stencil, 
the fifth  column  the standard diagonal BGK with  $\nu_0=0.05$, 
the sixth  column  MRT LBM scheme with the quartic condition (\ref{sigsig}) realized.
}
\end{centering} 
\end{table}

\begin{table}
\begin{centering}
\begin{tabular}{ccccccc}
$l$ & Bessel &     FD3-108 & FD5-108  &  BGK-108  &  LB-108 &  LB-108-q \cr
%%%  Bessel &     FD3-108 & FD5-108  &  BGK-108  &  LB-108   LB-108-q \cr
&&&Singlets&&&\cr
 1 &  14.68200  &  0.00165 &  -0.00052  &   0.00069  &  0.00070 &  0.00035 \cr
 2 &  49.21850  &  0.00628 &  -0.00104  &   0.00179  &  0.00189 &  0.00010 \cr
 3 & 103.49950  &  0.01382 &  -0.00175  &   0.00355  &  0.00377 & -0.00028 \cr
 4 & 177.52080  &  0.02399 &  -0.00230  &   0.00599  &  0.00640 & -0.00078 \cr
 5 & 271.28171  &  0.03665 &  -0.00244  &   0.00923  &  0.00989 & -0.00138 \cr
 6 & 384.78189  &  0.05198 &  -0.00202  &   0.01341  &  0.01442 & -0.00204 \cr
&&&Doublets&&&\cr
 1 &  26.37460 &   0.00410 &  -0.00027  &   0.00143 &  0.00147 &  0.00058 \cr
 2 &  40.70650 &   0.00662 &  -0.00029  &   0.00189 &  0.00197 &  0.00044 \cr
 3 &  57.58290 &   0.00943 &  -0.00053  &   0.00249 &  0.00262 &  0.00035 \cr
 4 &  76.93890 &   0.01251 &  -0.00066  &   0.00321 &  0.00339 &  0.00029 \cr
 5 &  98.72630 &   0.01601 &  -0.00086  &   0.00404 &  0.00428 &  0.00025 \cr
 6 & 122.90760 &   0.01979 &  -0.00101  &   0.00498 &  0.00530 &  0.00023 \cr
 7 & 149.45290 &   0.02380 &  -0.00115  &   0.00603 &  0.00643 &  0.00021  \cr
 8 & 178.33730 &   0.02814 &  -0.00121  &   0.00720 &  0.00768 &  0.00022  \cr
 9 & 209.54010 &   0.03281 &  -0.00126  &   0.00846 &  0.00905 &  0.00023  \cr
10 & 243.04340 &   0.03766 &  -0.00123  &   0.00983 &  0.01053 &  0.00022  \cr
11 & 278.83160 &   0.04274 &  -0.00118  &   0.01133 &  0.01215 &  0.00022  \cr 
 \end{tabular}
 \caption{Numerical study of the Stokes modes in a disk. 
 All simulations are performed with the same value $\nu_0 = 1 / \sqrt{108}$.
The second  column gives the theoretical values and the other the error 
between the numerical scheme and the theoretical value. 
The third column uses the present FD-LBM scheme with a three point
stencil for the evaluation of the gradients, 
the fourth column  the present FD-LBM scheme with a five point 
stencil, 
the fifth  column  the standard diagonal BGK with  $\nu_0 = 1 / \sqrt{108}$, 
the sixth  column  the MRT-LBM scheme with the quartic condition (\ref{sigsig}) not realized 
and the seventh   column the quartic version of the  MRT-LBM scheme 
when the condition (\ref{sigsig}) is realized. }
\end{centering} 
\end{table}

%%%%%%%%%%%%%%%%%%%%%%%%%%%%%%%%%%%%%%%%%%%%%%%%%%%%%%%%%%%%%%%  section 3-c 
\bigskip \bigskip   \noindent {\bf \large   3-c) $ \,\,\, $ Poiseuille flow }    
%%%%%%%%%%%%%%%%%%%%%%%%%%%%%%%%%%%%%%%%%%%%%%%%%%%%%%%%%%%%%%%%%%%%%%%%%%%%%%%%%%%%%%%  

\noindent 
Some simulations of Poiseuille flow have been performed to estimate the efficiency     %%%P
of the boundary conditions. We consider a channel with solid boundaries parallel
to the  $ \, Ox \, $ axis and periodic boundary conditions at the open ends. We adapt the 
boundary conditions to impose $v=0$ at $y_1=1-\xi$ and $y_2=N+\xi$. A uniform
body force parallel to $ \, Ox \, $  drives the flow. After enough time steps the stationary
flow is least square fit to a parabolic flow allowing to define ``experimental''
boundaries where the parabola goes to 0 at $y_{m1}=1-\xi_m$ and $y_{m2}=N+\xi_m$.
We show in Fig.~6 %%%  \ref{fig-poiseuille} 
the measured $\xi_m$ {\it vs}  the imposed $\xi$.

%
%%%%%%%%%%%%%%%%%%%%%%%%%%%%%%%%%     Figure 6        %%%%%%%%%%%%%%%%%%%%%%%%%%%%%% 
\begin{figure}[h!] \begin{center} 
\label{fig-poiseuille}  
\centerline { \includegraphics[width=.75 \textwidth] {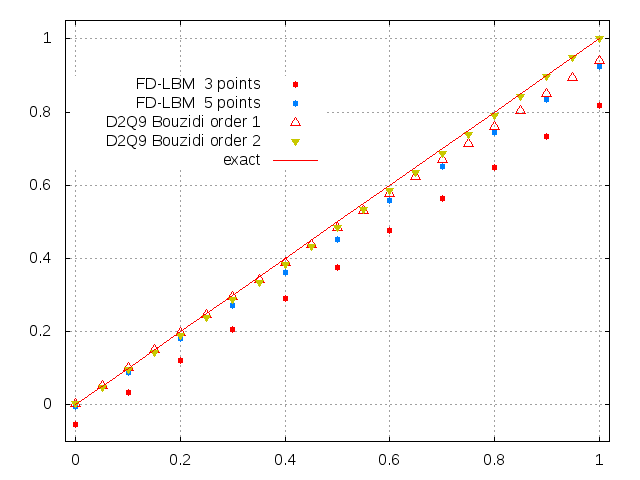} }
\caption{ Boundary conditions for a Poiseuille flow. 
Numerical location of the zero of velocity {\it{versus}} 
the imposed value $\xi$ of the boundary. 
In the transverse direction, the computational domain is composed by the interval $ \, [1-\xi ,\, 15+\xi] \,$
where $ \, \xi \,$ is the abscissa of the figure. The result is the location of the 
zero velocity value measured from a least square fit of the velocity profile. 
The result with the five points difference scheme is of good quality, 
comparable to what is obtained with the classical D2Q9 usual LBM scheme 
with first order extrapolation with the Bouzidi {\it et al.} algorithm. 
Observe that with a simple ``bounce-back'' boundary conditions, the result
would be a horizontal line at $ \, y = 0.5$. }
\end{center} \end{figure} 
%%%%%%%%%%%%%%%%%%%%%%%%%%%%%%%%%%%%%%%%%%%%%%%%%%%%%%%%%%%%%%%%%%%%%%%%%%%%%%%%%%%%
% 

%%%%%%%%%%%%%%%%%%%%%%%%%%%%%%%%%%%%%%%%%%%%%%%%%%%%%%%%%%%%%%%%%%  
\bigskip \bigskip   \noindent {\bf \large Conclusion }    
%%%%%%%%%%%%%%%%%%%%%%%%%%%%%%%%%%%%%%%%%%%%%%%%%%%%%%%%%%%%%%%%%% 

\noindent 
We have shown that the ACM proposal can be improved in two ways : reducing the velocity
dependence of the shear viscosity and diminishing the hyperviscosity with the use of
a stencil with more points. However when identical values of the long wave length
shear and bulk viscosities are chosen for the D2Q9 lattice Boltzmann model, the 
hyperviscosity is much smaller as can be seen in Fig.~5. %%  \ref{figd2q9}. 
An analogous analysis has also been performed for the three-dimensional model D3Q19.   %%%P

\noindent 
The present work needs to be complemented with detailed testing of situations where
nonlinear terms dominate to see the quality of simulations.
This will help decide how  many grid points in FD-LBM are needed to get 
comparable accuracy to what is given by  a  LBE-D2Q9 calculation.

%%%%%%%%%%%%%%%%%%%%%%%%%%%%%%%%%%%%%%%%%%%%%%%%%%%%%%%%%%%%%%%%%%%%%%%%%%%%%%%%%%%%%
\bigskip \bigskip   \newpage \noindent {\bf \large  References } 
%%%%%%%%%%%%%%%%%%%%%%%%%%%%%%%%%%%%%%%%%%%%%%%%%%%%%%%%%%%%%%%%%%%%%%%%%%%%%%%%%%%%%


\vspace{-.3cm} 
 
\begin{thebibliography}{99}

\bibitem{AOCR12}
P. Asinari, T. Ohwada, E. Chiavazzo, A.F. Di Rienzo. 
``Link-wise artificial compressibility method'', 
{\it Journal of Computational Physics}, vol.~\textbf {231},  p.~5109-5143, 2012. 

 \vspace{-.24cm}  
\bibitem{ADGL14}
A. Augier, F. Dubois, B. Graille and P. Lallemand. 
``On rotational invariance of Lattice Boltzmann schemes'', 
{\it Computers and Mathematics with Applications},  
vol.~\textbf {67},  p~239-255, 2014. 

 \vspace{-.24cm}  
\bibitem{BFL01} 
M. Bouzidi, M. Firdaous, P. Lallemand,  
``Momentum transfer of a Boltzmann-lattice fluid with boundaries'',
\textit {Physics of Fluids}, vol.~{\bf{13}}, p.~3452-3459, 2001.

 \vspace{-.24cm}  
\bibitem{CD98}
S. Chen, G. D. Doolen.
``Lattice Boltzmann method for fluid flows'', 
{\it Annual  Review of  Fluid Mechanics},  vol.~\textbf{30}, 
p.~329-364, 1998.

 \vspace{-.24cm}  
\bibitem{De02}  
P. J. Dellar. 
``Lattice Kinetic Schemes for Magnetohydrodynamics'', 
{\it Journal of Computational Physics},   vol.~\textbf {179},  p.~95-126, 2002.

 \vspace{-.24cm}  
\bibitem{De13}  
P. J. Dellar. 
``An interpretation and derivation of the lattice Boltzmann method using Strang splitting", 
{\it{Comput. Math. Applic.}}, vol. {\bf{65}}, p.~129-141, 2013.

 \vspace{-.24cm}  
\bibitem{Du08}
 F. Dubois.
``Equivalent  partial differential equations of a Boltzmann scheme'',
{\it Computers and mathematics with applications},  vol.~\textbf{55}, 
p.~1441-1449, 2008.

 \vspace{-.24cm}  
\bibitem{Du09}          
F. Dubois. 
``Third order equivalent equation of lattice Boltzmann scheme'', 
{\it Discrete and Continuous Dynamical Systems-Series A},  vol.~\textbf {23}, 
number 1/2,  p.~221-2482009, a special issue dedicated to Ta-Tsien Li 
on the occasion of his 70th birthday, doi: 10.3934/dcds.2009.23.221, 2009. 


 \vspace{-.24cm}  
\bibitem{DL09}           
F. Dubois, P. Lallemand. 
``Towards higher order lattice Boltzmann schemes'', 
{\it Journal of Statistical Mechanics: Theory and Experiment}, 
P06006 doi: 10.1088/1742-5468/2009/06/P06006, 2009. 


 \vspace{-.24cm}  
\bibitem{DL11}           
F. Dubois, P. Lallemand. 
``Quartic Parameters for Acoustic Applications of Lattice Boltzmann Scheme'', 
{\it Computers and mathematics with applications},  vol.~\textbf{61}, 
p.~3404-3416, 2011, doi:10.1016/j.camwa.2011.01.011.


 \vspace{-.24cm}  
\bibitem{DDH92}    
D. d'Humi\`eres. ``Generalized Lattice-Boltzmann Equations'',  in 
{\it Rarefied Gas Dynamics: Theory and Simulations}, 
vol.~\textbf {159} of \textit {AIAA Progress in
Astronautics and Astronautics}, p.~450-458, 1992.  

 \vspace{-.24cm}  
\bibitem{DGK}
D. d'Humi\`eres, I. Ginzburg, M. Krafczyk, P. Lallemand and L.S. Luo.
``Multiple-relaxation-time lattice Boltzmann models in three dimensions'', 
{\it Philosophical Transactions of the Royal Society},  London, vol.~\textbf{360}, 
p.~437-451, 2002.

 \vspace{-.24cm}  
\bibitem{JKL05}
M. Junk, A. Klar, L.S. Luo. ``Asymptotic analysis of the lattice Boltzmann equation'',
{\it Journal of Computational Physics}, vol.~\textbf {210},  p.~676-704, 2005. 

 \vspace{-.24cm}  
\bibitem{LL00}  
P. Lallemand, L-S. Luo. 
  ``Theory of the lattice Boltzmann method: 
   Dispersion, dissipation, isotropy, Galilean invariance, and stability'',
{\it Physical Review E}, vol.~{\bf 61}, p.~6546-6562, June 2000.  

 \vspace{-.24cm}  
\bibitem{OAKR14}
C. Obrecht, P. Asinari, F. Kuznik and J.J. Roux. 
``High-performance implementations and large-scale validation of the link-wise 
artificial compressibility method'', 
{\it Journal of Computational Physics}, vol.~{\bf 275}, p.~143-153, 2014.   

 \vspace{-.24cm}  
\bibitem{OA10}
T. Ohwada, P. Asinari.
``Artificial Compressibility Method Revisited: 
Asymptotic Numerical Method for Incompressible Navier-Stokes Equations'', 
{\it Journal of Computational Physics}, vol.~\textbf {229},  p.~1698-1723, 2010. 

\bibitem{OAY11}
T. Ohwada, P. Asinari, D. Yabusaki.
``Artificial Compressibility Method 
and Lattice Boltzmann Method: Similarities and Differences'',  
{\it Computers and Mathematics With Applications}, vol.~\textbf {61},  p.~3461-3474, 2011. 

 \vspace{-.24cm}  
\bibitem{WWLL13} %%   Jia Wang, Donghai Wang, Pierre Lallemand, Li-Shi Luo
J. Wang, D. Wang, P. Lallemand, L-S. Luo. 
 ``Lattice Boltzmann simulations of thermal convective flows in two dimensions'', 
{\it Computers and Mathematics with Applications}, 
  vol.~\textbf {65}, p.~262-286, 2013. 


 \vspace{-.24cm}  
\bibitem{Ye02}  
J. Yepez. 
``Quantum Lattice-Gas Model for the Burgers Equation'', 
{\it Journal of Statistical Physics}, 
vol.~\textbf {107},  p.~203-224, 2002.


\end{thebibliography}
\end{document}